\documentclass{gen-j-l}

\newtheorem{theorem}{Theorem}[section]

\newtheorem{conjecture}[theorem]{Conjecture}

\theoremstyle{definition}

\theoremstyle{remark}

\numberwithin{equation}{section}

\def\DJ{{\hbox{D\kern-.8em\raise.15ex\hbox{--}\kern.35em}}}

\def\NSERC{The author was supported in part by the
NSERC Grant A-5285.}

\def\bC{{\mbox{\bf C}}}
\def\bR{{\mbox{\bf R}}}

\begin{document}

\title[Atiyah's conjecture]
{Verification of Atiyah's conjecture for some nonplanar
configurations with dihedral symmetry}

\author[D.\v Z. \DJ okovi\'c]{Dragomir \v Z. \DJ okovi\'c}
\address{Department of Pure Mathematics, University of Waterloo,
Waterloo, Ontario, N2L 3G1, Canada}
\email{djokovic@uwaterloo.ca}
\thanks{\NSERC}

\date{}

\begin{abstract}
To an ordered $N$-tuple $(x_1,\ldots,x_N)$ of distinct points in
the three-dimensional Euclidean space, Atiyah has associated
an ordered $N$-tuple of complex homogeneous polynomials
$(p_1,\ldots,p_N)$ in two variables $x$ and $y$ of degree $N-1$,
each $p_i$ determined only up to a scalar factor.
He has conjectured that these polynomials are linearly independent.
In this note it is shown that Atiyah's conjecture is true
if $m$ of the points are on a line $L$ and the remaining
$n=N-m$ points are the vertices of a regular $n$-gon whose
plane is perpendicular to $L$ and whose centroid lies on $L$.
\end{abstract}

\maketitle



\date{}

\section{Introduction}

Let $(x_1,\ldots,x_N)$ be an ordered $N$-tuple of distinct
points in $\bR^3$. Each ordered pair
$(x_i,x_j)$ with $i\ne j$ determines a point
$$
\frac{x_j-x_i}{|x_j-x_i|}
$$
on the unit sphere $S^2\subset\bR^3$. Identify $S^2$ with the
complex projective line ${\bf CP}^1$ by using a stereographic
projection. We obtain a point 
$(u_{ij},v_{ij})\in{\bf CP}^1$ and a nonzero linear form
$l_{ij}=u_{ij}x+v_{ij}y\in\bC[x,y]$.
Define homogeneous polynomials $p_i\in\bC[x,y]$ of
degree $N-1$ by
\begin{equation} \label{poli}
p_i=\prod_{j\ne i}l_{ij}(x,y),\quad i=1,\ldots,N.
\end{equation}

\begin{conjecture} \label{conj-1} (Atiyah \cite{MA2})
The polynomials $p_1,\ldots,p_N$ are linearly independent.
\end{conjecture}

Atiyah \cite{MA1,MA2} observed that his conjecture is true
if the points $x_1,\ldots,x_N$ are collinear. He also
verified the conjecture for $N=3$.
Then Eastwood and Norbury \cite{EN} verified it for $N=4$.
In our previous note \cite{DZ} we verified this
conjecture for two special planar configurations of $N$ points.
For additional information on the conjecture (further
conjectures, generalizations, and numerical evidence)
see \cite{MA2,AS}.

Apart from the above mentioned result for arbitrary four points,
there are no results known for nonplanar configurations.
In this note we prove Atiyah's conjecture for the infinite
family of nonplanar configurations described in the abstract.

\section{Preliminaries}

Identify $\bR^3$ with $\bR\times\bC$ and denote the origin by $O$.
Following Eastwood and Norbury \cite{EN}, we make use of the Hopf map
$h:\bC^2\setminus\{O\}\to(\bR\times\bC)\setminus\{O\}$ defined by:
$$
h(z,w)=((|z|^2-|w|^2)/2,z\bar{w}).
$$
This map is surjective and its fibers are the circles
$\{(zu,wu):u\in S^1\}$, where $S^1$ is the unit circle in $\bC$.
If $h(z,w)=(a,v)$, we say that $(z,w)$ is a {\em lift} of $(a,v)$.

Let $x_i=(a_i,z_i)$.
For the sake of simplicity, we assume that if $i<j$ and
$z_i=z_j$ then $a_i<a_j$.
As the lift of the vector $x_j-x_i$, $i<j$, we choose
$$
\lambda_{ij}^{-1/2}\left(\lambda_{ij},\bar{z}_j-\bar{z}_i\right),
$$
where
$$
\lambda_{ij}=a_j-a_i+\sqrt{(a_j-a_i)^2+|z_j-z_i|^2}.
$$
Then
$$
\lambda_{ij}^{-1/2}\left(z_i-z_j,\lambda_{ij}\right),
$$
is a lift of $x_i-x_j$.
The corresponding linear forms are
\begin{eqnarray*}
l_{ij}(x,y) &=& \lambda_{ij}x+(\bar{z}_j-\bar{z}_i)y,\quad i<j; \\
l_{ij}(x,y) &=& (z_j-z_i)x+\lambda_{ji}y,\quad i>j.
\end{eqnarray*}

Define the binary forms $p_i$ by using (\ref{poli})
and the above expressions for the $l_{ij}$'s.
Atiyah's conjecture asserts that the
$N\times N$ coefficient matrix of these forms is nonsingular.

\section{Verification of the conjecture}

We shall prove Atiyah's conjecture for the configurations
of $N$ points satisfying the following two conditions:
\begin{itemize}
\item[(i)] The first $m$ points $x_1,\ldots,x_m$ lie on a line $L$.
\item[(ii)] The remaining $n=N-m$ points
$y_j=x_{m+j+1}$ $(j=0,1,\ldots,n-1)$
are the vertices of a regular $n$-gon
whose plane is perpendicular to $L$, and whose
centroid lies on $L$.
\end{itemize}

Without any loss of generality,
we may assume that $L=\bR\times\{0\}$
and that the $y_j$'s lie on the unit circle in $\{0\}\times\bC$.
Write $x_i=(a_i,0)$ for $i=1,\ldots,m$ and $y_j=(0,b_j)$
for $j=0,1,\ldots,n-1$.
We may also assume that $a_1<a_2<\cdots<a_m$
and that $b_j=-\zeta^j$, where $\zeta=e^{2\pi i/n}$.

The lifts of the nonzero vectors $x_j-x_i$,
$i,j\in\{1,\ldots,N\}$ are given in Table 1,
where we have set
$$
\lambda_i=a_i+\sqrt{1+a_i^2}.
$$

\vspace{3mm}

\begin{center}

{\bf Table 1: The lifts of the vectors $x_j-x_i$ }
$$
\begin{array}{l|c|r|c}
{\rm Vectors} & {\rm Index~ restrictions} &
{\rm Lifts}\quad\quad\quad & {\rm Linear~ forms} \\
\hline
x_r-x_i & i<r\le m & \left(2(a_r-a_i)\right)^{1/2}(1,0) & x \\
x_r-x_i & r<i\le m & \left(2(a_i-a_r)\right)^{1/2}(0,1) & y \\
y_s-y_j & s\ne j &
|b_s-b_j|^{1/2}\left(\frac{b_s-b_j}{|b_s-b_j|},1\right) &
\frac{b_s-b_j}{|b_s-b_j|}x+y \\
y_j-x_i & i\le m &
\lambda_i^{-1/2}(1,\lambda_i\bar{b}_j) & x+\lambda_i\bar{b}_jy \\
x_i-y_j & i\le m &
\lambda_i^{-1/2}(-\lambda_ib_j,1) & y-\lambda_ib_jx \\
\end{array}
$$
\end{center}
\vspace{3mm}

The associated polynomials $p_i$ (up to scalar factors)
are given by:
\begin{eqnarray*}
&& p_i(x,y)=x^{m-i}y^{i-1}\left(x^n-\lambda_i^ny^n\right),
\quad 1\le i\le m; \label{jed-m} \\
&& p_{m+j+1}(x,y)=\prod_{s\ne j}\left(x+
\frac{\bar{b}_s-\bar{b}_j}{|b_s-b_j|}y\right)\cdot
\prod_{i=1}^m(y-\lambda_ib_jx), \quad 0\le j<n. \label{jed-n}
\end{eqnarray*}

We now give the proof of our result.

\begin{theorem} \label{T-A}
Atiyah's conjecture is valid for configurations described above.
\end{theorem}

\begin{proof}
If $n=1$ or $2$, these configurations are planar and they
have been dealt with in \cite{DZ}. So, we assume that $n\ge3$.

Note that
$$
b_s-b_j=-2i\zeta^{j+s}\sin\frac{\pi(s-j)}{n}.
$$
After dehomogenizing the polynomials $p_i$ by setting $x=1$,
we obtain (up to scalar factors and ordering) the following polynomials:
\begin{eqnarray}
&& y^{i-1}(1-\lambda_i^ny^n),\quad 1\le i\le m; \label{pol-1} \\
&& f(\zeta^j y),\quad 0\le j<n, \label{pol-2}
\end{eqnarray}
where
\begin{equation} \label{pol-f}
f(y)=\prod_{s=1}^{n-1}(y-ie^{\pi is/n})\cdot\prod_{i=1}^m (y+\lambda_i).
\end{equation}

Denote by $\tilde{E}_k$ the $k$-th elementary symmetric function
of the $N-1$ numbers:
$$
\lambda_i,~(1\le i\le m);\quad -ie^{\pi is/n},~(1\le s\le n-1).
$$
By convention we set $\tilde{E}_0=1$ and $\tilde{E}_k=0$
if $k<0$ or $k\ge N$. Then
$$
f(y)=\sum_{k=0}^{N-1}\tilde{E}_{N-1-k}y^k.
$$
By factorizing $f$ over the real numbers, we see that all
coefficients of $f$ are positive.

Let $P$ be the coefficient matrix of the polynomials (\ref{pol-1}) and
(\ref{pol-2}). The top $m$ rows of $P$ form the submatrix
$$
\begin{pmatrix}
1 & 0 & 0 & \cdots & 0 & -\lambda_1^n & 0 & 0 & \ldots & 0 \cr
0 & 1 & 0 & \cdots & 0 & 0 & -\lambda_2^n & 0 & \ldots & 0 \cr
0 & 0 & 1 & \cdots & 0 & 0 & 0 & -\lambda_3^n & \ldots & 0 \cr
\vdots & & & & & & & & & \cr
\end{pmatrix}
$$
and the bottom $n$ rows the submatrix
$$
\begin{pmatrix}
\tilde{E}_{N-1} & \tilde{E}_{N-2} & \tilde{E}_{N-3} &
\cdots & \tilde{E}_1 & \tilde{E}_0 \cr
\tilde{E}_{N-1} & \tilde{E}_{N-2}\zeta & \tilde{E}_{N-3}\zeta^2 &
\cdots & \tilde{E}_1\zeta^{N-2} & \tilde{E}_0\zeta^{N-1} \cr
\tilde{E}_{N-1} & \tilde{E}_{N-2}\zeta^2 & \tilde{E}_{N-3}\zeta^4 &
\cdots & \tilde{E}_1\zeta^{2(N-2)} & \tilde{E}_0\zeta^{2(N-1)} \cr
\vdots &&&&& \cr
\end{pmatrix}.
$$

In order to compute $\det(P)$ we perform on $P$ successively
the following operations:
\begin{itemize}
\item[$\bullet$] Add the first column multiplied with
$\lambda_1^n$ to the $(n+1)$-st column.
\item[$\bullet$] Add the second column multiplied with
$\lambda_2^n$ to the $(n+2)$-nd column.
\item[$\vdots$]
\item[$\bullet$] Add the $m$-th column multiplied with
$\lambda_m^n$ to the $N$-th column.
\end{itemize}

By expanding the determinant of this new matrix
along the first $m$ rows, we obtain that
$$
|\det(P)|=c\prod_{k=0}^{n-1} f_k,
$$
where $c=n^{n/2}$ is the modulus of the determinant of the
matrix $(\zeta^{rs})$, $0\le r,s<n$, and
$$
f_k=\tilde{E}_k+\lambda_{m-k}^n\tilde{E}_{k+n}
+\lambda_{m-n-k}^n\lambda_{m-k}^n\tilde{E}_{k+2n}+\cdots,
\quad 0\le k<n.
$$
As the $\lambda_i$'s and the $\tilde{E}_k$'s are positive,
the proof is completed.
\end{proof}

\section{Comments on Atiyah and Sutcliffe conjecture}

Let us also state explicitly the stronger conjecture of
Atiyah and Sutcliffe \cite[Conjecture 2]{AS} for the
case of our configurations:
\begin{equation} \label{conj-2}
n^{n/2}\prod_{k=0}^{n-1} f_k\ge2^{\binom n 2 }
\prod_{i=1}^m (1+\lambda_i^2)^n,
\end{equation}
where, as in the proof above,
\begin{equation*}
f_k=\sum_{s\ge0}\left(\prod_{j=1}^s
\lambda_{N-jn-k}^n\right)\tilde{E}_{k+sn},\quad 0\le k<n.
\end{equation*}
Recall that $a_1<a_2<\cdots<a_m$ and, consequently,
$0<\lambda_1<\lambda_2<\cdots<\lambda_m$.

The substitution $(a_1,\ldots,a_m)\to(-a_m,\ldots,-a_1)$
corresponds to the reflection in the plane $\{0\}\times\bC$.
Consequently, the function
$$
\frac{\prod_{k=0}^{n-1}f_k}{\prod_{i=1}^m(1+\lambda_i^2)^n}
$$
is invariant under the transformation
$$
(\lambda_1,\ldots,\lambda_m)\to
(\lambda_m^{-1},\ldots,\lambda_1^{-1}).
$$

For $n=1$ the inequality (\ref{conj-2}) was proved in
\cite{DZ} in general, and for $n=2$ only in the special
(limit) case when all $\lambda_i$'s are equal.
One expects the inequality (\ref{conj-2}) to be strict
for all $n\ge3$ (see \cite[Section 4]{AS}).

Expand the two products in (\ref{pol-f}) separately:
\begin{eqnarray*}
&& \prod_{s=1}^{n-1}(y-ie^{\pi is/n})
=\sum_{j=0}^{n-1} c_j y^{n-1-j}, \\
&& \prod_{i=1}^m (y+\lambda_i)=\sum_{j=0}^m E_j y^{m-j}.
\end{eqnarray*}
The coefficients $c_j$, $0\le j<n$, and $E_j$, $0\le j\le m$, are
all positive. We also set $E_j=0$ if $j<0$ or $j>m$. Then
$$
\tilde{E}_k=\sum_{i=0}^{n-1}c_i E_{k-i}.
$$

In the limit case, when all $\lambda_k$'s are equal to
some $\lambda>0$, the inequality (\ref{conj-2}) specializes to
$$
n^{n/2}\prod_{k=0}^{n-1}\sum_{s\ge0}\lambda^{2sn+k}
\sum_{i=0}^{n-1}\binom{m}{k+sn-i} c_i\lambda^{-i}
\ge2^{\binom n 2}(1+\lambda^2)^{mn}.
$$
For $\lambda=0$ this gives
$$
n^{n/2}\prod_{k=0}^{n-1} c_k\ge2^{\binom n 2}.
$$

We conjecture that the following apparent strengthening
of (\ref{conj-2}) is valid:
\begin{equation} \label{spec}
\prod_{k=0}^{n-1}\frac{f_k}{c_k}
\ge\prod_{i=1}^{m}(1+\lambda_i^2)^{n}.
\end{equation}
When all $\lambda_k$'s are equal to
some $\lambda>0$, this becomes:
\begin{equation} \label{spec=}
\prod_{k=0}^{n-1}\sum_{s\ge0}\lambda^{2sn+k}
\sum_{i=0}^{n-1}\binom{m}{k+sn-i} c_i\lambda^{-i}
\ge \left( \prod_{k=0}^{n-1}c_k \right) \cdot (1+\lambda^2)^{mn}.
\end{equation}

If $n=3$ then $c_0=c_2=1$, $c_1=\sqrt{3}$ and
the inequality (\ref{spec}) takes the form:
\begin{equation} \label{n=3}
f_0f_1f_2\ge\sqrt{3}\prod_{i=1}^m(1+\lambda_i^2)^3,
\end{equation}
where
\begin{equation*}
f_k = \sum_{s\ge0}\left(\prod_{j=1}^s\lambda_{N-3j-k}^3\right)
(E_{3s+k}+\sqrt{3}E_{3s+k-1}+E_{3s+k-2}),
\end{equation*}
and (\ref{spec=}) the form:
$$
f_0f_1f_2\ge\sqrt{3}(1+\lambda^2)^{3m},
$$
where now
$$
f_k=\sum_{s\ge0}\lambda^{6s+k}\left[
\binom m {3s+k}+
\sqrt{3}\binom m {3s+k-1}\lambda^{-1}+
\binom m {3s+k-2}\lambda^{-2}\right].
$$

By using Maple, we have verified the last inequality for $m\le6$,
and, by using the invariance property mentioned above, it is easy to
verify (\ref{n=3}) for $m=2$.

\end{document}